\documentclass[12pt]{article}
\usepackage{mathtools}
\usepackage{ucs} 
\usepackage[utf8]{inputenc} 
\usepackage[T1,T2A]{fontenc}
\usepackage{graphicx}
\usepackage{wrapfig}
\usepackage{xcolor}
\usepackage{hyperref}
\usepackage{setspace}

 \title{\Large\bf Niels Henrik Abel and the Birth \\ of  Fractional Calculus}
 \author{Igor Podlubny$^1$,  Richard L. Magin$^2$, Iryna Trymorush$^1$}

\begin{document}
 
 \maketitle





 \begin{abstract}
 
\vbox{
\begin{flushright}
\textit{Dedicated to Professor Virginia Kiryakova} \\
\textit{on the occasion of her 65th birthday}\\
\textit{and the 20th anniversary of FCAA}
\end{flushright}
}

In his first paper on the generalization of the tautochrone problem, that was published in 1823, 
Niels  Henrik Abel presented a complete framework for fractional-order calculus,
and used the clear and appropriate notation for fractional-order integration and differentiation.

 \medskip

{\it MSC 2010\/}: Primary 26A33

 \smallskip

{\it Key Words and Phrases}:
fractional calculus, tautochrone, Caputo derivative, time scales, history of mathematics

 \end{abstract}

 \maketitle


\section{\uppercase{Introduction}}

\begin{wrapfigure}[20]{R}{6cm}
\begin{center}
      \includegraphics[width=6cm]{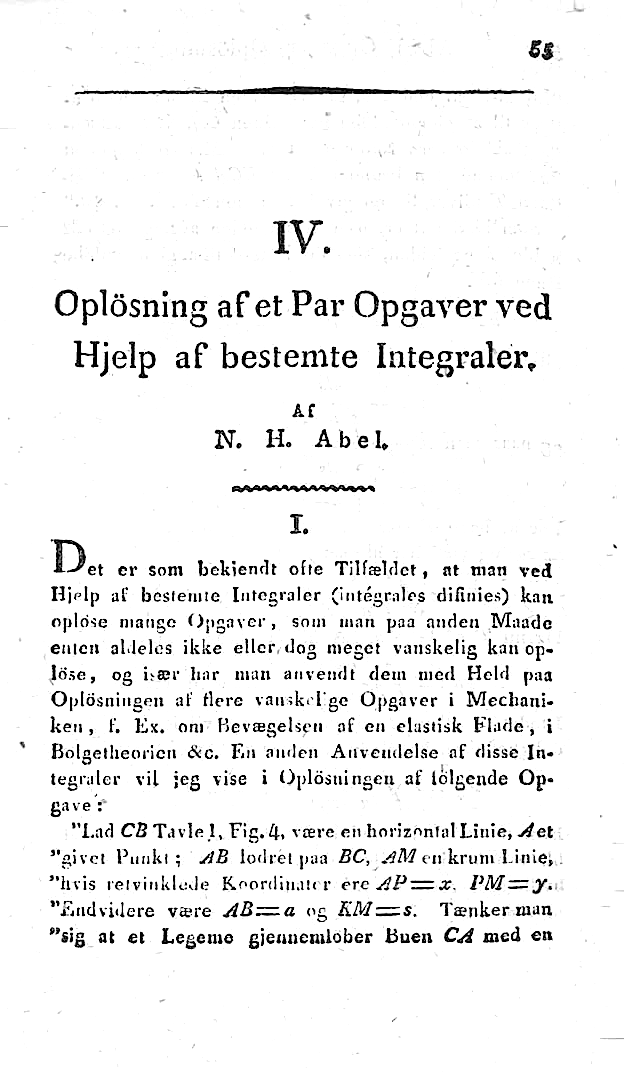}
\end{center}
 \end{wrapfigure}

Niels Henrik Abel's life was too short. He was born on August 5, 1802 (that is, exactly 215 years ago) and he passed away on April 6, 1829, at the age of less than 27 years. Nowadays, many graduate students of a similar age are still pursuing  their degrees. 
Nevertheless, during his short career he had made several monumental mathematical contributions. 

His most famous result is, of course, the proof that a general algebraic equations of the fifth degree cannot be solved in radicals. 

In addition, in every textbook on calculus we can find the Abel's tests for the convergence of an infinite series (for number series, and for power series). 

Another of his contributions was the doubly-periodic functions that are now called the Abel elliptic functions. 

His name is also associated with several smart transformations, such as Abel's summation by parts, as well as with Abel's inequality,  Abel's binomial expansion, the Abelian means for divergent series, and  many other elegant results and formulas, without which one cannot imagine many fields of the today's mathematics.

However, one of Abel's most fascinating inventions, 
the formulas for fractional-order integration and differentiation,
has remained practically unnoticed. 
We found these formulas  working on a totally different topic, 
where we needed to take a closer look at the tautochrone problem. 
In reading Abel's papers on this topic we discovered that in solving the generalization of the tautochrone problem 
Niels Henrik Abel had also developed a complete framework of what is now called the fractional calculus, or differentiation and integration of arbitrary real order. 

In his first article (\cite{Abel-1823}, 1823), the first page of which is shown here,   
N.~H.~Abel introduced 
fractional-order integration in the form that is currently known 
as the Riemann-Liouville fractional integral, 
and fractional-order differentiation in the form that is 
currently known as the Caputo fractional derivative. 
Abel's solution of the considered problem is, in fact, 
the proof that these two operators are mutually inverse. 
This means that N.~H.~Abel, who was only 21 years of age at the time of  the publication of his paper,
 was the father  of  the complete fractional-order calculus framework.

For  reasons that remained unknown, in his second paper (\cite{Abel-1826}, 1826)
Abel abandoned this line of theoretical development and 
its applications, and never returned to it. 
In the following we take a look at Abel's reasoning in more detail.

\section{\uppercase{Abel's 1823 paper}}

~
\vspace*{-3ex}

\begin{wrapfigure}[6]{R}{0.2\textwidth}
\begin{center}
      \includegraphics[width=0.22\textwidth]{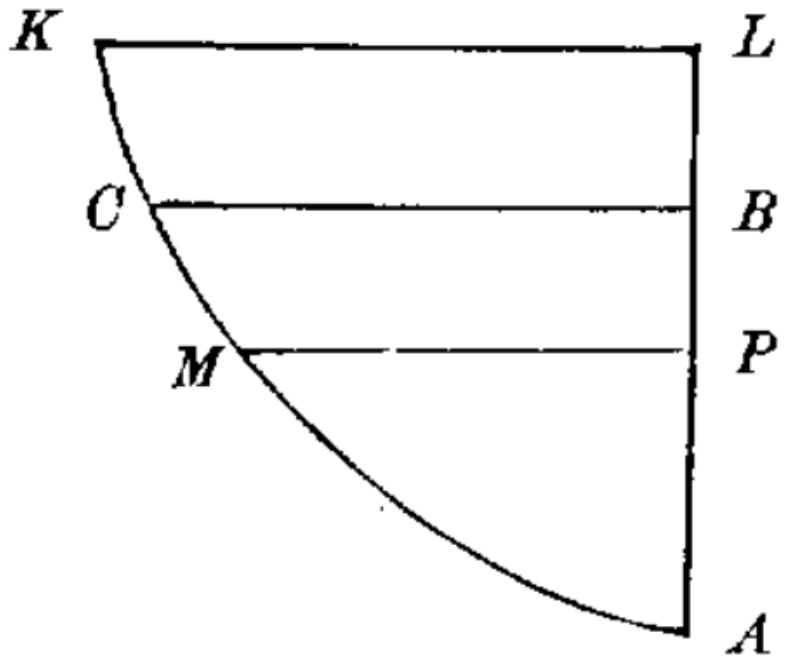}
\end{center}
\end{wrapfigure}

\noindent
First of all, it should be noted that Abel considered a more general problem than the tautochrone 
(in the quotation below, which is translated and transcribed from \cite{Abel-1823}, 
 $\psi a$ means in today's notation $\psi (a)$):

\smallskip

\begin{quotation}
``Suppose that $CB$ is a horizontal line, $A$~is a setpoint, $AB$ is perpendicular to $BC$, $AM$ ia a curve with rectangular coordinates $AP=x$, $PM=y$. Moreover $AB=a$, $AM=s$. It is known that as a  body moves along an arc $CA$, when the initial velocity is zero, that the time $T$, which is necessary for the passage, depends on the shape of the curve, and on $a$. One has to find the definition of a curve $KCA$, for which the time $T$ is a given function of $a$, for example $\psi a$.''
\end{quotation}

\smallskip

\noindent
The picture on the right appeared in the French translation in \cite{Oeuvres-1881}. \\
For this problem Abel obtains the following equation:
$$
	\psi a = \int \frac{ ds }{ \sqrt{ a - x  } } 
	\quad 
	(\mbox{from } x = 0, \mbox{ to } x = a),
$$
and then continues (we added the equation numbers for convenience): 

\smallskip

\begin{quotation}
``Instead of solving this equation, I will show how one can derive $s$ from the more general equation
\begin{equation}\label{eq:Abel-equation-n}
	\psi a = \int \frac{ ds }{ ( a - x )^n }
	\quad 
	(\mbox{from } x = 0, \mbox{ to } x = a)
\end{equation}
where $n$ has to be less than 1 to prevent the infinite integral between two limits; $\psi a$ is an arbitrary function that is not infinite, when $a$ equals to zero.''
\end{quotation}

\smallskip

\noindent
Abel looks for the unknown function $s(x)$ in the form of a power series, and after term-by-term operations with that series uses the properties of the gamma function (paying credits to A.-M. Legendre~\cite{Legendre}, who studied and summarized the properties of the Euler's gamma function, introduced the notation $\Gamma(z)$ and gave the function its name). After several pages of manipulations with the power series he arrives at 
the solution of equation~(\ref{eq:Abel-equation-n}):
\begin{equation}
	s = \frac{ \sin n \pi }{ \pi } x^n 
	\int \frac{ \psi( xt ) dt }{ ( 1 - t )^{ 1 - n } }
	\quad 
	( t = 0, \mbox{ } t = 1),
\end{equation}
and he calls this ``m\ae rkv\ae rdige Theorem'' --  a remarkable theorem.

Indeed, it is  remarkable, since Abel further examines the obtained formulas from other viewpoints and mentions that 
``\dots $s$ can be expressed in the different way'':
\begin{equation}
	s = \frac{1}{ \Gamma( 1 - n ) } \int^n \psi x . dx^n = \frac{1}{ \Gamma( 1 - n ) } \frac { d^{ - n } \psi x }{ dx^{ - n } }.
\end{equation}
In spite of the archaic historical notation, we easily recognize that here Abel uses two expressions for the fractional-order integral: one is the derivative of negative order, and the other is the symbol that later appears in Liouville's works, namely $\int^n \psi x . dx^n$. This means that Abel understood that he unified the notions of integration and differentiation, and that he extended them to non-integer orders!

Writing the equation (\ref{eq:Abel-equation-n}) as 
\begin{equation}\label{eq:Abel-equation-n-today}
	\psi (t) = 
	\int\limits_{ 0 }^{  t }
	 \frac{ s'(x) dx }{ ( t - x )^n },
\end{equation}
we observe that Abel's equation is nothing else but (up to a constant coefficient) 
the nowadays famous and widely used Caputo fractional derivative of real order $n$  ($0 \leq n < 1$),
and the solution of  Abel's equation is simply the inverse operation -- 
fractional-order integral of the same order $n$. 

In his 1823 paper \cite{Abel-1823} Abel freely uses non-integer orders even in the text, for example: 
``Differentieres V\ae rdien for $s$ $n$Gange, saa faaer man\ldots'' is literally translated as
``En differentiat $n$ fois de suite la valeur de $s$, on obtient...'' \cite{Oeuvres-1881}, or
``Differentiating $n$ times the value of $s$, one obtains...'', and $n$ here is non-integer!
Moreover, returning to the consideration of the original equation, Abel writes:

\smallskip

\begin{quotation}
``If $n = \frac{1}{2}$, one obtains
$$
	\psi a = \int \frac{ ds }{\sqrt{a - x}  },
		\quad 
	(x = 0, \, x = a),
$$
and
$$
	s = \frac{1}{\sqrt{\pi}} \frac{ d^{ - \frac{1}{2} } \psi x }{ dx^{ - \frac{1}{2} } }  
	= \frac{1}{\sqrt{\pi}} \int^{ \frac{1}{2} } \psi x . dx^{ \frac{1}{2} }.
$$
This is the equation of the desired curve, when the time~is~$\, \psi a$.''
\end{quotation}

\smallskip

\noindent
And at this moment Abel is  ready to simply invert the fractional integral by applying  fractional differentiation to both sides of the equation and by using the understanding (and the fact), that fractional-order integration and differentiation are mutually inverse operations:

\smallskip

\begin{quotation}
``One obtains from this equation
$$
	\psi x = \sqrt{\pi} \frac{ d^{ \frac{1}{2} } s }{ dx^{ \frac{1}{2} } },
$$
therefore: if the curve equation is $s = \psi x$, then the time that the body uses to pass the arc, 
the height of which is $a$, equals
$\displaystyle{\sqrt{\pi} \frac{ d^{\frac{1}{2}} \psi a  }{ da^{\frac{1}{2} } }}$.''
\end{quotation}

\section{\uppercase{Abel's 1826 paper}}

~
\vspace*{-3ex}

\begin{wrapfigure}[8]{R}{0pt}
      \includegraphics[width=0.25\textwidth]{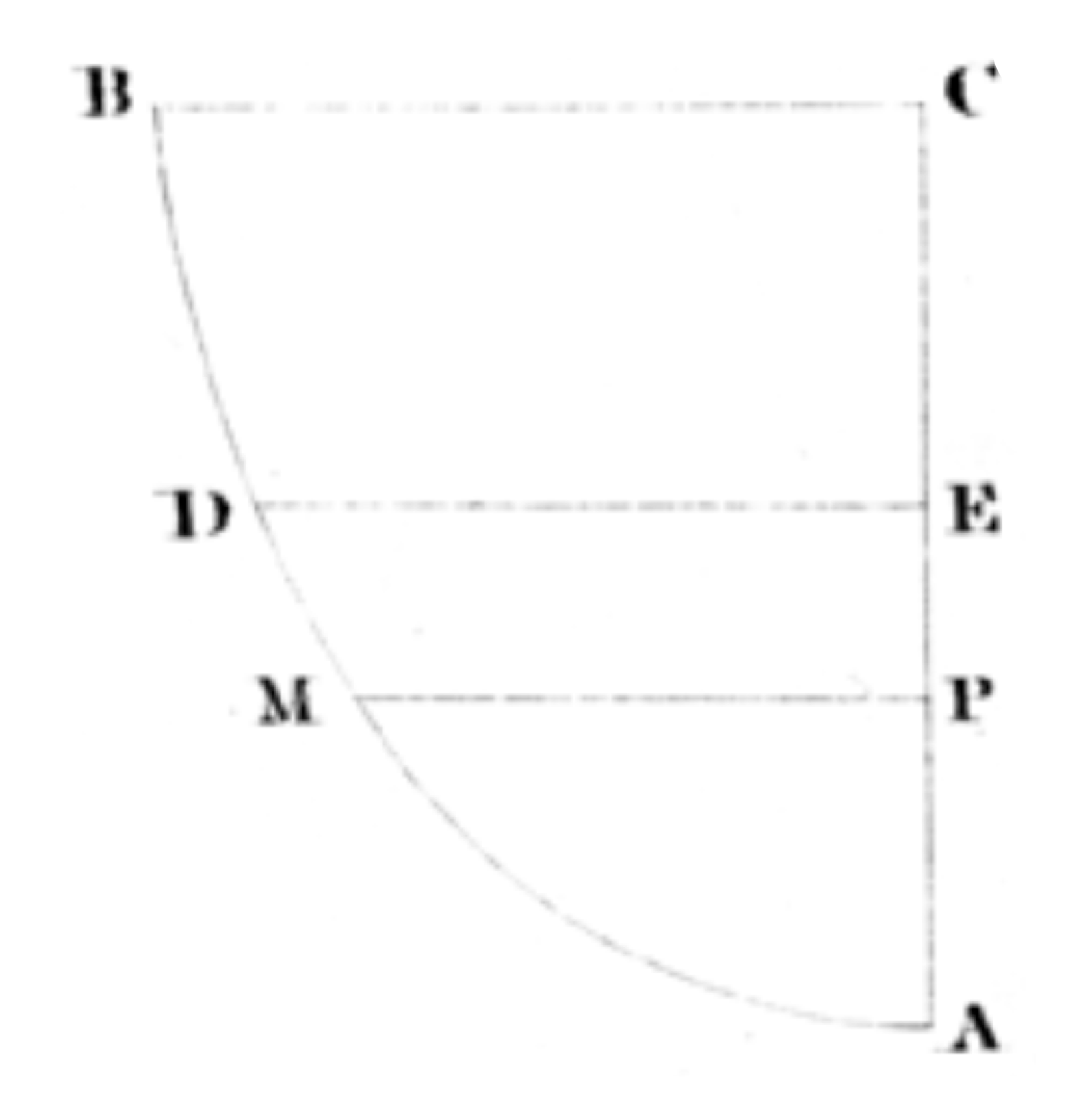}
\end{wrapfigure}

\noindent
We speculate that Abel was not very satisfied with the manipulations of the power series in his 1823 paper \cite{Abel-1823}; maybe he considered them too lengthy, or not well justified, or maybe simply not so elegant. In either case, in his paper dated 1826 \cite{Abel-1826}, that is three years after his first paper on the tautochrone problem, he presented a different method of solution, which was based on using the properties of the Euler's gamma function. The figure on the right is from the original paper.

First, he derives the following formula (written here in Abel's notation):
\begin{equation} \label{eq:Abel-FO-transform}
	f \, x = \frac{\sin n\pi}{\pi}  \int_{0}^{x}\frac{da}{(x-a)^{1-a}}  
		\int_{0}^{a}\frac{f'z \,dz}{(a-z)^n},
\end{equation}
and today we easily recognize here the fractional-order integral of the Caputo fractional-order derivative.  Then Abel uses this formula for obtaining the solution 
of equation (\ref{eq:Abel-equation-n}): 

\smallskip

\begin{quotation}
``Multiplying this equation [\textit{that is, (\ref{eq:Abel-equation-n})}] by 
$\displaystyle{\frac {\sin n\pi}{\pi} \cdot \frac{da}{(x-a)^{1-n}}}$, 
and integrating it from $a = 0$ to $a = x$, one obtains
$$
	\frac{\sin n\pi}{\pi} \int_{0}^{x}\frac{\varphi a . da}{(x-a)^{1-n}} =
	\frac{\sin n\pi}{\pi} 
	\int_{0}^{x} \frac{da}{(x-a)^{1-n}} 
	\int_{0}^{a} \frac{ds}{(a-x)^{n}}
$$ 
Therefore, accordingly to the equation (1) [\textit{here (\ref{eq:Abel-FO-transform})}]
$$
	s = 	\frac{\sin n\pi}{\pi} 
		\int_{0}^{x}\frac{\varphi a  \, da}{(x-a)^{1-n}}. \mbox{''}
$$
\end{quotation}

\smallskip

\noindent
In this paper Abel also considers the case of $n=\frac{1}{2}$, as he did three years earlier:

\smallskip

\begin{quotation}
``Considering for now $n = \frac{1}{2}$, one obtains 
$$
	\varphi \, a = \int_{0}^{a}\frac{ds}{ {\sqrt{}} \, (a - x) }
$$ 
\noindent
and
$$
	s = \frac{1}{\pi} \int_{0}^{x} \frac{\varphi a \, da}{ {\sqrt{}} \, ( x - a)  }.  \mbox{''}
$$
\end{quotation}

\smallskip

\noindent
The 1826 paper \cite{Abel-1826} is shorter, and the derivation of the solution 
of  equation (\ref{eq:Abel-equation-n}) is concise and looks more elegant;
but the fractional-order calculus, which was, in fact, 
developed in Abel's 1823 paper, disappeared\ldots 
\linebreak
Abel shortened and improved the solution of a particular problem, 
but lost the whole theory of the fractional-order calculus.

\section{\uppercase{Space and time scales}}

The tautochrone problem, considered by Abel, also has another aspect, 
that remained unnoticed. In fact, it provides an example of a whole family of dynamical processes, which, so to say, `live' in the same time scale. Indeed, the time necessary to reach the destination from each point  of the arc of particular shape is the same. 
However, the spatial scale is different for  each starting point, as the distance 
passed from the start to the end is different for each point on the arc.

\section{\uppercase{Conclusion}}

It is not clear why Niels Henrik Abel abandoned the direction of research so nicely 
formed in his 1823 paper \cite{Abel-1823}, and one can only guess the reasons. 
Abel had all the elements of the fractional-order calculus there: 
the idea of fractional-order integration and differentiation, 
the mutually inverse relationship between them, 
the understanding that fractional-order differentiation and integration 
can be considered as the same generalized operation,
and even the unified notation for differentiation and integration of arbitrary real order. 

It is easy to propose the following hypothesis regarding why the mathematical community did not 
pay proper attention to Abel's invention of the fractional-order calculus. 
The original paper was published in the Danish language in the Norwegian journal 
\emph{Magazin for Naturvidenskaberne}, so the readership was very limited. 
Thus, Abel's 1823 paper was not included in his posthumous ``\OE uvres compl\`etes de N. H. Abel'' 
compiled by B. Holmboe and published in Christiania in 1839 \cite{Oeuvres-1839}, 
while a translation of Abel's 1826 paper \cite{Abel-1826} was included there as Chapter IV.  
Later, in 1881, L. Sylow and S. Lie compiled another 
``\OE uvres compl\`etes de Niels Henrik Abel'' \cite{Oeuvres-1881}, 
also printed in Christiania; they included a translation of Abel's 1823 \hbox{paper} as Chapter II,
and a translation of Abel's 1826 paper as Chapter IX. 
However, apparently the initial picture of Abel's scientific work 
created by Holmboe's collection \cite{Oeuvres-1839} lasted too long,
and Abel's 1823 paper remained in the shadow of his elegant 1826 paper that became 
known to the mathematical community earlier,  
as it was published in the already famous Crelle's journal.

In Abel's 1823 paper \cite{Abel-1823} we find the notation that was later used by Liouville for fractional-order integration, and the definition of fractional differentiation that is now called the Caputo fractional differentiation. It took 144 years until Professor Michele Caputo came up with his definition in his paper in 1967 \cite{Caputo-1967} and in his book \cite{Caputo-1969} in 1969.  
	However, as Hans \hbox{Selye} wrote, the discovery  
\emph{``is not to see something first, but to establish solid connections 
between the previously known and the hitherto unknown''}~\cite{Selye-Stress-Book}.
Abel's discovery of the fractional calculus remained unnoticed in his days (even by himself), so Professor Caputo's name is therefore properly associated with this definition, as in his works he demonstrated numerous  applications of fractional-order differentiation to viscoelasticity, geophysics, and other fields.



\section*{Acknowledgments}


\noindent
This work was  supported in parts by grants VEGA 1/0908/15,  
APVV-14-0892, SK-PL-2015-0038, ARO W911NF-15-1-0228.


 \bigskip \smallskip

{ \it

 \noindent
$^1$ 
Institute of Control and Informatization \\ 
of Production Processes, BERG Faculty, \\
Technical University of Kosice, \\ 
B. Nemcovej 3, 04200 Kosice, Slovakia \\ 
e-mail: igor.podlubny@tuke.sk, iryna.trymorush@tuke.sk\\[2ex]
$^2$ 
Department of Bioengineering \\ 
University of Illinois at Chicago, \\
851 S Morgan St, SEO 218 \\
Chicago, IL 60607, USA \\
e-mail: rmagin@uic.edu\\[12pt]
}

\newpage
\section*{Appendix I}

\begin{center}

{\large\bf
N. H. Abel\\
Solution of some problems using definite integrals. \\
Aargang I, Bind 2, Christiania 1823. }

\bigskip

{Translation by: \\
Iryna Trymorush\\
Technical University of Kosice, Slovakia\\
email: iryna.trymorush@tuke.sk}
\date{}

\end{center}


\begin{abstract}

This is the supplement \#1 to the paper: 
Podlubny, I., Magin, R. L., Trymorush I., 
``Niels Henrik Abel and the birth of fractional calculus'',
\emph{Fractional Calculus and Applied Analysis}, vol.~20, no.~5, pp.~1068--1075, 2017
(\url{https://doi.org/10.1515/fca-2017-0057}). 

Hereby a translation to English of the first part of the important Abel's paper 
(Abel, N. H., 
Opl\"osning af et par opgaver ved hjelp af bestemte integraler,
\emph{Magazin for Naturvidenskaberne}, Aargang I, Bind 2, Christiania, 1823)
is provided. The second part of the original Abel's paper has been omitted, 
since it is not related to the historical development of fractional-order calculus.

To preserve the original flavor and the historical development of mathematical notation,
the formulas are typeset in the original  manner.
\end{abstract}

\newpage

\begin{wrapfigure}[5]{R}{0pt}
\includegraphics[width=0.15\textwidth]{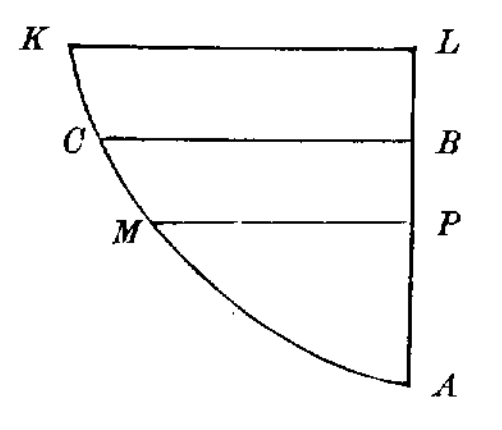}
\end{wrapfigure}

It is well known that many problems are solved by definite integrals  which otherwise can not be solved, or at least are very difficult to solve. They have especially been applied with advantage to the solution of several difficult problems of mechanics, for example, to the motion of an elastic surface, problems of wave theory, etc. I will show a new application by solving the following problem.\\
Suppose, that $CB$ is a horizontal line, $A$ is a setpoint, $AB$ is perpendicular to $BC$, $AM$ is a curve with rectangular coordinates $AP=x$, $PM=y$. Moreover $AB=a$, $AM=s$. Let us suppose that the body moves along an arc $CA$, the initial velocity is zero, the time $T$, which is necessary for the passage, depending on the shape of the curve, and on $a$. We are talking about the definition of a curve $KCA$, for which the time $T$ is a given function of $a$, for example $\psi a$.\\
Assume that $h$ is the velocity of a point mass at $M$, and $t$ is a time required for it to travel along the arc $CM$, it is known that
$$
	h = \sqrt{ BP } = \sqrt{ ( a - x ) }
$$
$$
	dt = - \frac{ ds }{ h },
$$
therefore
$$
	dt = - \frac{ ds }{ \sqrt{ ( a - x ) } },
$$
and integrating
$$
	t = - \int \frac{ ds }{ \sqrt{ ( a - x ) } }.
$$
To find $T$ it is needed to take the integral from $x=a$ to $x=0$, thus
$$
	T = \int_{ x = 0 }^{ x = a } \frac{ ds }{ \sqrt{ ( a - x ) } }.
$$
If we define $T$ as $\psi a$, the equation becomes
$$
	\psi a = \int_{ x = 0 }^{ x = a } \frac{ ds }{ \sqrt{ ( a - x ) } }.
$$
Instead of solving this equation, I will show how one can derive $s$ from the more general equation
$$
	\psi a = \int_{ x = 0 }^{ x = a } \frac{ ds }{ ( a - x )^n },
$$
where $n$ has to be less than 1 to prevent the infinite integral between two limits; $\psi a$ is an arbitrary function that is not infinite, when $a$ equals to zero.\\
Let us say
$$
	s = \sum a^{ ( m ) } x^m,
$$
where $\sum a^{ ( m ) } x^m$ has following value:
$$
	\sum a^{ ( m ) } x^m = a^{ ( m' ) } x^{ m' } + a^{ ( m'') } x^{ m'' } + a^{ ( m''' ) } x^{ m'''} + ... .
$$
Differentiating one obtains
$$
	ds = \sum m a^{ ( m ) } x^{ m - 1 } dx,
$$
so
$$
	\frac{ ds }{ ( a - x )^n} = \frac{ \sum ma^{ ( m ) } x^{ m - 1 } dx }{ ( a - x )^n } = \sum ma^{ ( m ) } \frac{ x^{ m - 1 } dx }{ ( a - x )^n }.
$$
Integrating
$$
	\int_{ x = 0 }^{ x = a } \frac{ ds }{ ( a - x )^n } = \int_{ x = 0 }^{ x = a } \sum ma^{ ( m ) } \frac{ x^{ m - 1 } dx }{ ( a - x )^n }.
$$
Since
$$
	\int \sum ma^{ ( m ) } \frac{ x^{ m - 1 } dx }{ ( a - x )^n } = \sum  ma^{ ( m ) } \int \frac{ x^{ m - 1 } dx }{ ( a - x )^n },
$$
so as
$\int_{ x = 0 }^{ x = a } \frac{ ds }{ ( a - x )^n } = \psi a$:
$$
	\psi a = \sum ma^{ ( m ) } \int_0^a \frac { x^{ m - 1 } dx }{ ( a - x )^n }.
$$
The value of the integral
$$
	\int_0^a \frac{ x^{ m - 1 } dx }{ ( a - x )^n }
$$
is easy to find in the following way: if put $x=at$, will get
$$
	x^m = a^m t^m,
	mx^{ m - 1 } dx = ma^m t^{ m - 1 } dt,
$$
$$
	( a - x )^n = ( a - at )^n = a^n ( 1 - t )^n,
$$
so
$$
	\frac{ mx^{ m - 1 } dx }{ ( a - x )^n } = \frac{ ma^{ m - n } t^{ m - 1 } dt }{ ( 1 - t )^n },
$$
integrating
$$
	m \int_0^a \frac{ x^{ m - 1 } dx }{ ( a - x )^n } = ma^{ m - n } \int_0^1 \frac{ t^{ m - 1 } dt }{ ( 1 - t )^n }.
$$
Note that
$$
	\int_0^1 \frac{ t^{ m - 1 } dt }{ ( 1 - t )^n } = \frac{ \Gamma( 1 - n ) \Gamma m}{\Gamma( m - n + 1 )},
$$
where $\Gamma m$ is a function defined by the equations
$$
	\Gamma( m + 1 ) = m \Gamma m, 
	\Gamma( 1 ) = 1 . 
	\footnote{
	The properties of this remarkable function have been largely developed by M. Legendre in his book, "Exercises of integral calculus T.1 and T.2".}
$$
Now, substituting this value for the integral $\int_0^1 \frac{ t^{ m - 1 } dt }{ ( 1 - t )^n }$, and noticing that $m \Gamma m = \Gamma( m + 1 )$ one obtains
$$
	m \int_0^a \frac{ x^{ m - 1 } dx }{ ( a - x )^n } = \frac{ \Gamma( 1 - n ) \Gamma( m + 1 ) }{ \Gamma( m - n + 1 ) } a^{ m - n }.
$$
Substituting this value in the expression for $\psi a$, one obtains
$$
	\psi a = \Gamma( 1 - n ) \sum a^{ ( m ) }  a^{ m - n } \frac{ \Gamma( m + 1 ) }{ \Gamma( m - n + 1 ) }.
$$
Let us consider
$$
	\psi a = \sum \beta^{ ( k ) } a^{  k },
$$
one obtains
$$
	\sum \beta^{ ( k ) } a^{ k } = \sum \frac{ \Gamma( 1 - n ) \Gamma( m + 1 ) }{ \Gamma( m - n + 1 ) } a^{ ( m ) } a^{ m - n }.
$$

For this equation to be satisfied $ m - n = k$, therefore $m = n + k$, and 
$$
	\beta^{ ( k ) } = \frac{ \Gamma( 1 - n ) \Gamma( m + 1 ) }{ \Gamma( m - n + 1 ) } a^{ ( m ) } = \frac {\Gamma( 1 - n ) \Gamma( n + k + 1 ) }{ \Gamma( k + 1 ) } a^{ ( m ) },
$$
so
$$
	a^{ ( m ) } = \frac{ \Gamma( k + 1 ) }{ \Gamma( 1 - n ) \Gamma( n + k + 1 ) } \beta^{ ( k ) }.
$$
Then
$$
	\int_0^1 \frac{t^k dt}{ ( 1 - t )^{ 1 - n } } = \frac{ \Gamma n . \Gamma( k + 1 ) }{ \Gamma( n + k + 1 ) },
$$
so
$$
	a^{ ( m ) } = \frac{ \beta^{ ( k ) } }{ \Gamma n. \Gamma( 1 - n ) } \int_0^1 \frac{ t^k dt}{ ( 1 - t )^{ 1 - n } }.
$$
By multiplying by $x^m = x^{ n + k }$ one obtains
$$
	a^{ ( m ) } x^m = \frac{ x^n }{ \Gamma n. \Gamma( 1 - n ) } \int_0^1 \frac{ \sum \beta^{ ( k ) } (xt)^k dt}{ ( 1 - t )^{ 1 - n } },
$$
so that
$$
	\sum a^{ ( m ) } x^m = \frac{ x^n }{ \Gamma n. \Gamma( 1 - n ) } \int_0^1 \frac{ \sum \beta^{ ( k ) } (xt)^k dt}{ ( 1 - t  )^{ 1 - n } }.
$$
But $\sum a^{ ( m ) } x^m = s$, $\sum \beta^{ ( k ) } (xt)^k = \psi( xt )$, so
$$
	s = \frac{ x^n }{ \Gamma n. \Gamma( 1 - n ) }\int_0^1 \frac{ \psi( xt ) dt }{ ( 1 - t )^{ 1 - n }}.
$$
Then, noting that $\Gamma n. \Gamma( 1 - n ) = \frac{ \pi }{ \sin n \pi }$, wiil find
$$
	s = \frac{ \sin n \pi . x^n }{ \pi } \int_0^1 \frac{ \psi( xt ) dt }{ ( 1 - t )^{ 1 - n } }.
$$
From which follows the remarkable theorem:\\
if
$$
	\psi a = \int_{ x = 0 }^{ x = a } \frac{ ds }{ ( a - x )^n },
$$
then one obtains
$$
	s = \frac{ \sin n \pi }{ \pi } x^n \int_0^1 \frac{ \psi( xt ) dt }{ ( 1 - t )^{ 1 - n } }.
$$
Use this theorem for the equation
$$
	\psi a = \int_{ x = 0 }^{ x = a } \frac{ ds }{ \sqrt{ ( a - x ) } }.
$$
In this case  $n = \frac{1}{2}$, so $ 1 - n = \frac{1}{2}$ and therefore
$$
	s = \frac{ \sqrt x }{ \pi } \int_0^1 \frac{ \psi (xt) dt}{ \sqrt{ 1 - t } }.
$$
This equation defines the arc $s$ of the curve, as a function of the corresponding abscissa $x$; we can easily obtain an equation between rectangular coordinates by noting that  $ds^2 = dx^2 + dy^2$.
\\
\\
Let us apply the previous solution to some special cases.
\begin{itemize}
\item 
Find a curve with the property that the time that the body expends to move along an arbitrary arc is proportional to the $n^{th}$ power of the height traveled by body.\\
In such case $\psi a = c a^n$, where $c$ is a constant value, \\
so $\psi (xt) = cx^nt^n$, as the result we get:
$$
	s = \frac{ \sqrt x }{ \pi } \int_0^1 \frac{ cx^n t^n dt }{ \sqrt{ 1 - t } } 
	   = x^{ n + \frac{1}{2} } \frac{ c } { \pi } \int_0^1 \frac{ t^n dt }{ \sqrt{ 1 -  t} },
$$
thus using
$$
	\frac{c}{\pi} \int_0^1 \frac{ t^n dt }{ \sqrt{ 1 - t } } = C,
$$
one obtains
$$
	s = Cx^{  n + \frac{1}{2} };
$$
since
$$
	ds = ( n + \frac{1}{2} ) Cx^{  n - \frac{1}{2} } dx,
$$
and
$$
	ds^2 = ( n + \frac{1}{2} )^2 C^2 x^{ 2n - 1 } dx^2
	         = dy^2 + dx^2,
$$
from this deduce for $( n + \frac{1}{2} )^2 C^2 = k$
$$
	dy = dx \sqrt{ kx^{ 2n - 1 } - 1 }.
$$
The equation of the desired curve becomes
$$
	y = \int dx \sqrt{ kx^{ 2n - 1 } - 1 }.
$$
If I set $n=\frac{1}{2}$, and having $x^{2n-1} = 1$, one obtains
$$
	y = \int dx \sqrt{ k - 1 } = k' + x \sqrt{ k - 1 }.
$$
Consequently, the desired curve is a straight line.

\newpage

\item
Find the isochron equation.\\
Since the time has to be independent of the distance traveled, having $\psi a = c$ hence one obtains
$$
	s = \frac{ \sqrt{ x } }{ \pi } c \int_0^1 \frac{ dt }{ \sqrt{ 1 - t } },
$$
therefore
$$	
	s = k \sqrt{x},
$$
where
$$
	k = \frac{ c }{ \pi } \int_0^1 \frac{ dt }{ \sqrt{ 1 - t } }.
$$
This is the known cycloid equation.\\
It is known that if
$$
	\psi a = \int_{ x = 0 }^{ x = a } \frac{ ds }{ ( a - x )^n },
$$
and
$$
	s = \frac{ \sin n \pi }{ \pi } x^n \int_0^1 \frac{ \psi(xt) dt }{ ( 1 - t )^{ 1 - n } }.
$$
Since $s$ can be expressed in a different way, because of the singularity, it is known that
$$
	s = \frac{1}{ \Gamma( 1 - n ) } \int^n \psi x . dx^n = \frac{1}{ \Gamma( 1 - n ) } \frac { d^{ - n } \psi x }{ dx^{ - n } },
$$
so if
$$
	\psi a = \int_{ x = 0 }^{ x = a } ds( a - x )^n,
$$
and
$$
	s = \frac{1}{ \Gamma( 1 + n ) } \frac{ d^n \psi x}{ dx^n };
$$
in other words
$$
	\psi a = \frac{1}{ \Gamma( 1 + n ) } \int_{ x = 0 }^{ x =  a } \frac{ d^{ n + 1 } \psi x }{ dx^{ n + 1 } } ( a - x )^n dx.
$$
This assertion is easily proved in the following way. If we put
$$
	\psi x = \sum a^{ ( m ) } x^{ m }
$$
and differentiating one obtains
$$
	\frac{ d^n \psi x }{ dx^n } = \sum a^{ ( m ) } m ( m - 1 )( m - 2 )...( m - k + 1 ) x^{ n - k };
$$
but
$$
	m ( m - 1 )( m - 2 )...( m - k + 1 ) = \frac{\Gamma( m + 1 ) }{ \Gamma( m - k + 1 ) },
$$
so
$$
	\frac{ d^k \psi x }{ dx^k } = \sum a^{ ( m ) } \frac{ \Gamma( m + 1 ) }{ \Gamma( m - k + 1 ) } x^{ m - k }.
$$
Now
$$
	\frac{ \Gamma( m + 1 ) }{ \Gamma( m - k + 1 ) } = \frac{1}{ \Gamma( - k ) } \int_0^1 \frac{ t^m dt }{ ( 1 - t  )^{ 1 + k } },
$$
therefore
$$
	\frac{ d^k \psi x }{ dx^k } = \frac{1}{ x^k \Gamma( - k ) } \int_0^1 \frac{\sum a^{ ( m ) } (xt)^m dt }{ ( 1 - t )^{ 1 + k } };
$$
but $\sum a^{ ( m ) } (xt)^m = \psi( xt )$, thus
$$
	\frac{ d^k \psi x }{ dx^k } = \frac{1}{ x^k \Gamma( - k ) } \int_0^1 \frac{ \psi (xt) dt }{ ( 1 - t )^{ 1 + k } }.
$$
By putting $k = - n$, one obtains
$$
	\frac{ d^{ - n } \psi x }{ dx^{ - n } } = \frac{ x^n }{ \Gamma n } \int_0^1 \frac{ \psi (xt) dt }{ ( 1 - t )^{ 1 - n } }.
$$
It has been seen that
$$
	s = \frac{ x^n }{ \Gamma n . \Gamma( 1 - n ) } \int_0^1 \frac{ \psi( xt ) dt }{ ( 1 - t )^{ 1 - n } },
$$
so
$$
	s = \frac{1}{ \Gamma( 1 - n ) } \frac{ d^{ - n }  \psi x }{ dx^{ - n } },
$$
if
$$
	\psi a = \int_{ x = 0 }^{ x = a } \frac{ ds }{ ( a - x )^n },	
$$
Q.E.D.\\
Differentiatin $n$ times the value of $s$, one obtains
$$
	\frac{ d^n s }{ dx^n } = \frac{1}{ \Gamma( 1 - n ) } \psi x,
$$
and therefore define $s=\varphi x$,
$$
	\frac{ d^n \varphi a }{ da^n } = \frac{1}{ \Gamma( 1 - n ) } \int_0^a \frac{ \varphi ' x . dx }{ ( a - x ) ^n }.
$$
It should be noted that in the previous case $n$ must always be less than one.\\
If $n = \frac{1}{2}$, one obtains
$$
	\psi a = \int_{ x = 0 }^{ x = a } \frac{ ds }{\sqrt {( a - x )} }
$$
and
$$
	s = \frac{1}{\sqrt \pi} \frac{ d^{ - \frac{1}{2} } \psi x }{ dx^{ - \frac{1}{2} } }  = \frac{1}{\sqrt \pi} \int^{ \frac{1}{2} } \psi x . dx^{ \frac{1}{2} }.
$$
This is the equation of the desired curve, when the time is $\psi a$.\\
One obtains from this equation
$$
	\psi x = \sqrt{\pi} \frac{ d^{ \frac{1}{2} } s }{ dx^{ \frac{1}{2} } },
$$
therefore:\\
if the curve equation is $s = \varphi x$, the time that the body uses to travel along the arc, the height of which is $a$, equals
$\sqrt{\pi} \frac{ d^{ \frac{1}{2}} \varphi a }{ da^{\frac{1}{2} } }$.\\

Finally, I note that in the same way, accordingly to the equation
$$
	\psi a = \int_{ x = 0 }^{ x = a } \frac{ ds }{ ( a - x )^n }.
$$
I found $s$, also from the equation
$$
	\psi a = \int \psi(xa) fx . dx.
$$
I found the functions $\phi$, if $\varphi$ and $f$ are given functions, and the integral can be taken between any limits; but the solution of this problem is too long to be given here.
\end{itemize}

Value of expression  $\varphi ( x + y \sqrt{ - 1 } ) + \varphi( x - y \sqrt{ - 1 } )$.\\
When $\varphi$ is an algebraic, logarithmic, exponential or circular function, as it is known, we can always express the real value
$\varphi( x + y \sqrt{ - 1 } ) + \varphi( x - y \sqrt{ - 1 } )$ 
in real and finite form.
If from other side, $\varphi$ keeps its generality, we have not, as far as I know, been able to express it in real and finite form.
It can be done with the help of integrals defined as follows.\\
If we deduce from Taylor's theorem $\varphi (x + y \sqrt{ - 1 } )$ and $\varphi( x - y  \sqrt { - 1 } )$, one obtains
$$
	\varphi( x + y \sqrt{ - 1 } ) = 
						\varphi x +
						\varphi ' x.y\sqrt{ - 1 } - 
						\frac{ \varphi'' x }{ 1.2 } y^2 - 
						\frac{ \varphi''' x }{ 1.2.3 } y^3 \sqrt{ - 1 } + 
						\frac{ \varphi'''' x }{ 1.2.3.4 } y^4...
$$
$$
	\varphi( x - y \sqrt{ - 1 } ) = 
						\varphi x -
						\varphi ' x.y\sqrt{ - 1 } - 
						\frac{ \varphi'' x }{ 1.2 } y^2 + 
						\frac{ \varphi''' x }{ 1.2.3 } y^3 \sqrt{ - 1 } + 
						\frac{ \varphi'''' x }{ 1.2.3.4 } y^4...
$$

\newpage
\section*{Appendix II}

\begin{center}

{\large\bf
N. H. Abel\\
Solution of a mechanical problem.\\
Journal f{\"u}r die Reine und Angewandte Mathematik,  Band I, 153--157, Berlin, 1826.}

\bigskip

{Translation by: \\
Iryna Trymorush\\
Technical University of Kosice, Slovakia\\
email: iryna.trymorush@tuke.sk}\date{}

\end{center}

\begin{abstract}

This is the supplement \#2 to the paper: 
Podlubny, I., Magin, R. L., Trymorush I., 
``Niels Henrik Abel and the birth of fractional calculus'',
\emph{Fractional Calculus and Applied Analysis}, vol.~20, no.~5, pp.~1068--1075, 2017
(\url{https://doi.org/10.1515/fca-2017-0057}). 

Hereby a translation to English of the important Abel's paper 
(Abel, N. H.,  Aufl{\"o}sung einer mechanischen ausgabe.
\emph{Journal f{\"u}r die Reine und Angewandte Mathematik},  
Band I, 153--157, Berlin, 1826) is provided.

To preserve the original flavor and the historical development of mathematical notation,
the formulas are typeset in the original \linebreak manner.
\end{abstract}

\newpage

\begin{wrapfigure}[5]{R}{0pt}
\includegraphics[width=0.15\textwidth]{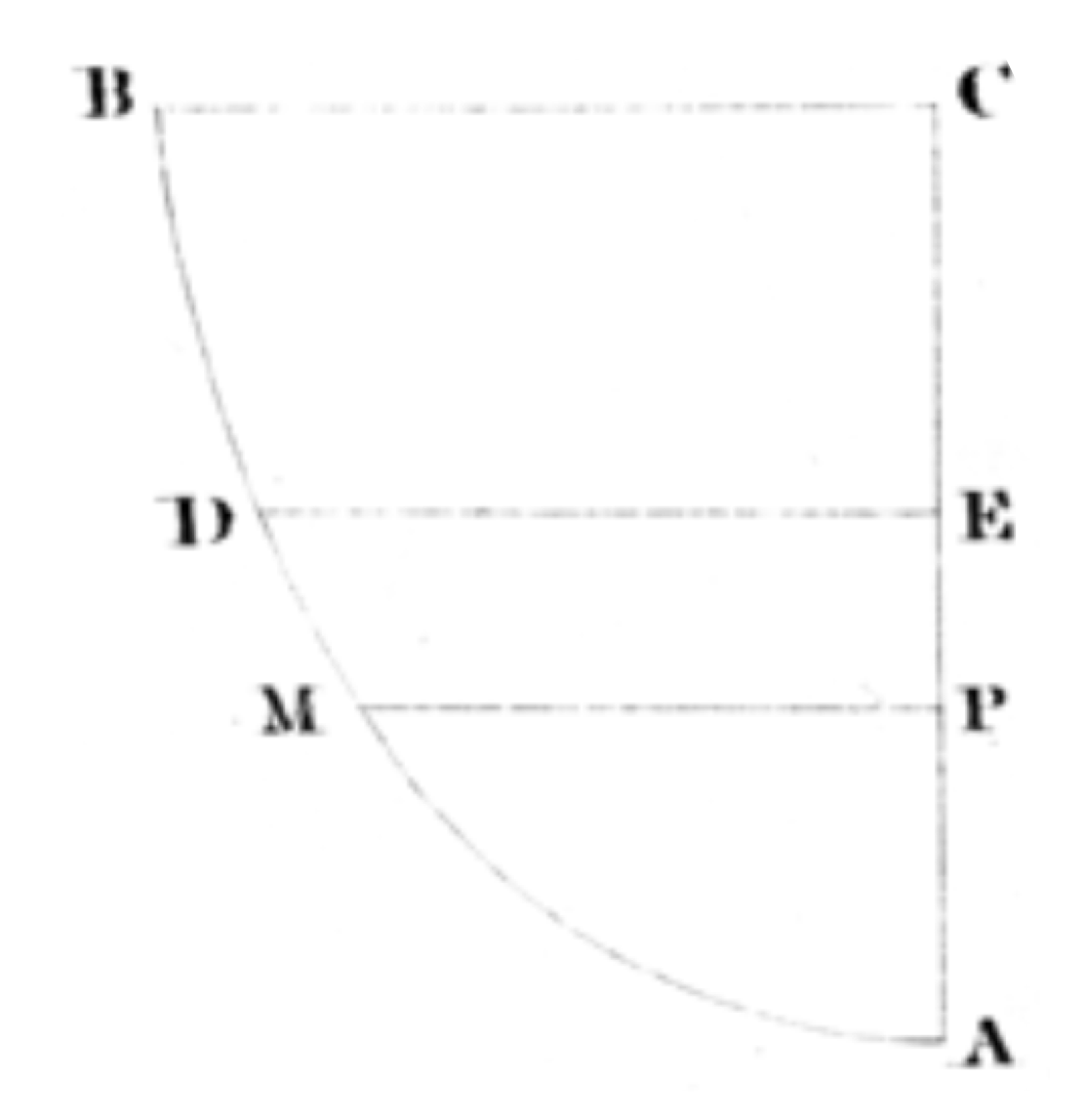}
\end{wrapfigure}

\noindent
Let $BDMA$ be an arbitrary curve.
The line $BC$  is a horizontal and $CA$ vertical.
Let us suppose that a point set in motion by the action of gravity moves on this curve, with an arbitrary point $D$ as its starting point.
Let $\tau$ be a time, which ends, when the moving body achives the point $A$, and let  $a$  be the height  $EA$.
The time $\tau$ is some function of $a$ and depends on the curve's form.
Conversely, the shape of the curve depends on this function.
We will examine how, with the aid of a defined integral, we can find the equation of the curve, for which $\tau$  is a continuous function of $a$.\\
Let $AM=s$, $AP=x$, and $t$ be the time that the body employs to traverse the arc $DM$. 
According to the rules of mechanics\\
$- \frac{ds}{dt} = \sqrt{ ( a - x ) }$, therefore $dt = - \frac{ds}{ \sqrt{ ( a - x ) } }$.\\
So, if we take an integral from  $x = a$ to  $x = 0$;\\
\begin{displaymath}
	\tau = - \int_{a}^{0} \frac{ds}{\sqrt{ ( a - x ) } }
	       = +\int_{0}^{a}\frac {ds}{\sqrt{ ( a - x ) } },
\end{displaymath}
\noindent
where $\int_{\alpha}^{\beta}$
means, that the limits of the integral are $x = \alpha$ and $x = \beta$.
Let us suppose for now, that
$$
	\tau = \varphi(a)
$$ 
is a given function, where
$$
	\varphi(a) = \int_{0}^{a} \frac{ds}{ \sqrt{ ( a - x ) } }
$$
\noindent
is the equation that can help find $s$ in $x$.
Instead of this equation we will consider the more general one 
$$
	\varphi(a) = \int_{0}^{a} \frac{ds}{ ( a - x )^n }
$$ 
to find $s$ in $x$.
Let us define the gamma function $\Gamma( \alpha )$ by
$$
	\Gamma( \alpha ) = \int_{0}^{1} dx ( \log \frac{1}{x} )^{ \alpha - 1 }, 
$$
which is known to solve
$$
	\int_{0}^{1} y^{ \alpha - 1 } ( 1 - y )^{ \beta - 1} dy = \frac{ \Gamma( \alpha ) . \Gamma( \beta ) }{ \Gamma( \alpha + \beta ) },
$$
where $\alpha$ and $\beta$ must be greater than zero.
Let $\beta = 1 - n$, then one obtains
$$
	\int_{0}^{1} \frac{ y^{ \alpha - 1 } dy }{ ( 1 - y )^{n} } =	\frac{ \Gamma( \alpha ) . \Gamma( 1 - n ) }{ \Gamma( \alpha + 1 - n ) },
$$
and defining $z = a y$, one obtains
$$
	\int_{0}^{a} \frac{ z^{ \alpha - 1 } dz }{ ( a - z )^{n} } = \frac{ \Gamma( \alpha ) . \Gamma( 1- n ) }{ \Gamma( \alpha + 1 - n ) } a^{ \alpha - n }.
$$
Multiplying by $\frac{da}{ ( x - a )^{ 1 - n } }$ and integrating from $a = 0$ to $a = x$, one obtains:
$$
	\int_{0}^{x} \frac{da}{ ( x - a )^{ 1 - n } } . \int_{0}^{a} \frac{ z^{  \alpha - 1 } dz}{ ( a - z )^n } =
	\frac{ \Gamma( \alpha ) . \Gamma( 1 - n ) }{ \Gamma( \alpha + 1 - n ) } . 
	\int_{0}^{x} \frac{ a^{ \alpha - n } da }{ ( x - a )^{ 1 - n } }.
$$
\noindent
If we now define $a=xy$, one obtains 
$$
	\int_{0}^{x} \frac{ a^{  \alpha - n } da}{ ( x - a )^{ 1 - n } } =
	x^\alpha \int_{0}^{1} \frac{ y^{ \alpha - n } dy }{ ( 1 - y )^{ 1 - n } } =
	x^\alpha \frac{ \Gamma( \alpha + 1 - n ) . \Gamma(n) }{ \Gamma( \alpha + 1 ) },
$$
\noindent
therefore
$$
	\int_{0}^{x} \frac{da}{ ( x - a )^{ 1 - n } } . 
	\int_{0}^{a}\frac{ z^{ \alpha - 1 } dz}{ ( a - z )^n }=
	\Gamma(n) . \Gamma( 1 - n ) . \frac{ \Gamma(\alpha) }{ \Gamma( \alpha + 1 ) } x^\alpha.
$$
\noindent
However, according to a known property of the gamma function, $\Gamma(\alpha)$ 
$$
	\Gamma(\alpha + 1)=\alpha \Gamma(\alpha);
$$
\noindent
therefore substituting I obtain:
$$
	\int_{0}^{x} \frac{da}{ ( x - a )^{ 1 - n } } . 
	\int_{0}^{a} \frac{ z^{ \alpha - 1 } dz }{ ( a - z )^n } =
	\frac{ x^\alpha }{\alpha} \Gamma(n) . \Gamma( 1 - n ).
$$
Multiplying by $\alpha . \varphi(\alpha) . d\alpha$  and integrating with respect to $\alpha$, one finds:
$$
	\int_{0}^{x} \frac{da}{ ( x - a )^{ 1 - n } } . 
	\int_{0}^{a} \frac{ (\int \varphi(\alpha) .  \alpha z^{ \alpha - 1 } d\alpha) dz}{ ( a - z )^n } = \Gamma(n) . \Gamma( 1 - n ) \int\varphi \alpha . x^\alpha d\alpha.
$$
Let
$$
	\int\varphi (\alpha)  x^\alpha d\alpha = f(x), 
$$
\noindent
then on differentiating
$$
	\int\varphi (\alpha) . \alpha . x^{\alpha-1} d\alpha = f'(x),
$$
$$
	\int\varphi (\alpha) . \alpha .  z^{\alpha-1} d\alpha = f'(z);
$$
\noindent
so
$$
	\int_{0}^{x} \frac{da}{ ( x - a )^{ 1 - n } } . 
	\int_{0}^{a} \frac{ f'(z) dz}{ ( a - z )^n } = \Gamma(n) . \Gamma(1-n) . f(x),
$$
\noindent
or since
$$
	\Gamma(n) . \Gamma(1-n) = \frac{\pi}{\sin n\pi},
$$ 
\begin{equation}
	f(x) = \frac{\sin n\pi}{\pi}  \int_{0}^{x}\frac{da}{(x-a)^{1-n}} . 
		\int_{0}^{a}\frac{f'(z)dz}{(a-z)^n}.
\end{equation}
Using this equation it will be easy to derive the value of $s$ from
$$
	\varphi(a) = \int_{0}^{a} \frac{ds}{ ( a - s )^n }.
$$
Multiplying this equation by $\frac {\sin n\pi}{\pi} . \frac{da}{(x-a)^{1-n}}$, 
and integrating it from $a = 0$ to $a = x$, one obtains
$$
	\frac{\sin n\pi}{\pi} \int_{0}^{x}\frac{\varphi (a) da}{(x-a)^{1-n}} =
	\frac{\sin n\pi}{\pi} 
	\int_{0}^{x} \frac{da}{(x-a)^{1-n}} 
	\int_{0}^{a} \frac{ds}{(a-x)^{n}},
$$ 
\noindent
therefore accordingly to equation (1) 
$$
	s = 	\frac{\sin n\pi}{\pi} 
		\int_{0}^{x}\frac{\varphi(a) da}{(x-a)^{1-n}}.
$$
\noindent
Let us consider for now $n = \frac{1}{2}$, so one obtains 
$$
	\varphi(a) = \int_{0}^{a}\frac{ds}{ \sqrt( a - x ) },
$$ 
\noindent
and
$$
	s = \frac{1}{\pi} \int_{0}^{x} \frac{\varphi(a) da}{ \sqrt( x - a ) }.
$$
\noindent
This equation gives the arc $s$ as a function of the abscissa $x$, and consequently the curve is entirely determined.
We will apply the expression found to some examples.
\begin{itemize}
\item [I.]
Let
$$
	\varphi(a) = \alpha_0 a^{\mu_0} + \alpha_1 a^{\mu_1} + \ldots + \alpha_m a^{\mu_m}
		      = \sum(\alpha a^\mu),
$$
and the value $s$ will be
$$
	s=\frac{1}{\pi}\int_{0}^{x}\frac{da}{\sqrt(x-a)} . \sum(\alpha a^\mu) =
		\frac{1}{\pi}
		\sum \left(
				\alpha \int_{0}^{x} \frac{a^\mu da}{\sqrt(x-a)}
			\right).
$$
\noindent
Defined as $a=xy$, one obtains 
$$
	\int_{0}^{x} \frac{a^\mu da}{\sqrt(x-a)} =
	x^{\mu+\frac{1}{2}} \int_{0}^{1} \frac{y^\mu dy}{\sqrt(1-y)} =
	x^{\mu+\frac{1}{2}} . \frac{\Gamma(\mu+1) . \Gamma(\frac{1}{2})}
					      {\Gamma(\mu+\frac{3}{2})},
$$
\noindent
therefore
$$
	s =	\frac{\Gamma(\frac{1}{2})}{\pi}
		\sum 
		\frac{\alpha \Gamma(\mu+1)}{\Gamma(\mu+\frac{3}{2})}
		x^{\mu+\frac{1}{2}},
$$
\noindent
and using $\Gamma(\frac{1}{2}) = \sqrt{\pi}$
$$
	s = \sqrt{ \left( \frac{x}{\pi} \right) } 
	\left[
		\alpha_0 \frac{\Gamma(\mu_0+1)}{\Gamma(\mu_0+\frac{3}{2})} x^{\mu_0} + 
		\alpha_1 \frac{\Gamma(\mu_1+1)}{\Gamma(\mu_1+\frac{3}{2})} x^{\mu_1} + 
		\ldots + 
		\alpha_m \frac{\Gamma(\mu_m+1)}{\Gamma(\mu_m+\frac{3}{2})} x^{\mu_m}
	\right].
$$
If we assume, for example, that $m = 0$, $\mu_0 = 0$, that is to say, that the curve sought is isochronous,
$$
	s = 
		\sqrt{\left( \frac{x}{\pi} \right)} . \alpha_0 \frac{\Gamma(1)}{\Gamma(\frac{3}{2})} = 
		\frac{\alpha_0}{\frac{1}{2}\Gamma(\frac{1}{2})} \sqrt{\left( \frac{x}{\pi} \right)}  =
		\frac{2 \alpha_0}{\pi}\sqrt{x},
$$ 
i.e. $s=\frac{2 \alpha_0}{\pi}\sqrt x$ is the known equation of the cycloid.
\item[II.]
Let's consider 
$$
\begin{array}{ll}
\varphi a 	\mbox{ from } a=0 \mbox{ to } a=a_0, 	& \mbox{equals  } \varphi_0 a\\
\varphi a 	\mbox{ from } a=a_0 \mbox{ to } a=a_1, 	& \mbox{equals  } \varphi_1 a\\
\varphi a 	\mbox{ from } a=a_1 \mbox{ to } a=a_2, 	& \mbox{equals  } \varphi_2 a\\
\ldots \ldots \ldots \ldots \ldots \ldots \ldots \ldots &  \\
\varphi a 	\mbox{ from } a=a_{m-1} \mbox{ to } a=a_m, & \mbox{equals  } \varphi_m a,\\
 \end{array}
 $$
\noindent
then one obtains \\
$\pi s=\int_{0}^{x} \frac{\varphi_0 a  . da}{\sqrt (a-x)}$, from $x=0$ to $x=a_0$\\
$\pi s=\int_{0}^{a_0} \frac{\varphi_0 a . da}{\sqrt (a-x)}
	+\int_{a_0}^{x} \frac{\varphi_1 a . da}{\sqrt (a-x)}$, from $x=a_0$ to $x=a_1$\\
$\pi s=\int_{0}^{a_0} \frac{\varphi_0 a . da}{\sqrt (a-x)}
	+\int_{a_0}^{a_1} \frac{\varphi_1 a . da}{\sqrt (a-x)}
	+\int_{a_1}^{x} \frac{\varphi_2 a . da}{\sqrt (a-x)}$, from $x=a_1$ to $x=a_2$\\
\ldots \ldots \ldots \ldots \ldots \ldots \ldots \ldots 
\ldots \ldots \ldots \ldots \ldots \ldots \ldots \ldots  \\	
$\pi s=\int_{0}^{a_0} \frac{\varphi_0 a . da}{\sqrt (a-x)}+\int_{a_0}^{a_1} \frac{\varphi_1 a . da}{\sqrt (a-x)}+...+\int_{a_{m-2}}^{a_{m-1}} \frac{\varphi_{m-1} a . da}{\sqrt (a-x)} + \int_{a_{m-1}}^{x} \frac{\varphi_{m} a . da}{\sqrt (a-x)}$, from $x=a_{m-1}$, to $x=a_m$,\\
if we notice that the functions $\varphi_0 a$, $\varphi_1 a$, $\varphi_2 a$, \ldots, $\varphi_m a$ 
have to be such that 
$$
	\varphi_0(a_0) = \varphi_1(a_0), 
	\varphi_1(a_1) = \varphi_2(a_1), 
	\varphi_2(a_2)=\varphi_3(a_2),
$$
then the functions $\varphi a$ have to be continuous.
\end{itemize}

\end{document}